\documentclass[11pt]{amsart}
\usepackage[T1]{fontenc}
\usepackage{lmodern,microtype,amsmath,amssymb,amsthm}
\usepackage{graphicx}
\graphicspath{{figures/}}
\usepackage[margin=1.15in]{geometry}
\usepackage[hidelinks]{hyperref}
\hypersetup{pdftitle={The four-dimensional Borsuk number is five or six},pdfauthor={Logan R. Chalmers}}
\newtheorem{theorem}{Theorem}
\newtheorem{proposition}[theorem]{Proposition}
\newtheorem{lemma}[theorem]{Lemma}
\DeclareMathOperator{\diam}{diam}
\DeclareMathOperator{\tr}{tr}
\DeclareMathOperator{\conv}{conv}
\newcommand{\R}{\mathbb R}
\newcommand{\ip}[2]{\langle #1,#2\rangle}
\newcommand{\norm}[1]{\lVert #1\rVert}
\newcommand{\clB}{\overline B}
\title{The four-dimensional Borsuk number is five or six}
\author{Logan R. Chalmers}
\address{University of Otago, Dunedin 9016, New Zealand}
\email{logan.chalmers@postgrad.otago.ac.nz}
\date{}
\begin{document}
\begin{abstract}
Can every set of diameter one in four dimensions be split into at most
five parts of smaller diameter?
Three parts suffice in the plane and four in three-dimensional space.
In four dimensions five parts are sometimes necessary. We prove that
six always suffice, improving the previous upper bound of seven.
For small circumradius we use six cones. For large
circumradius five balls meet every diametral pair of a compact set.
At intermediate radii we bound the diameter of each part by the diameter
of a rational polytope. The theorem is formalised in Lean.
\end{abstract}
\maketitle

\section{Introduction}
Can every set of diameter one in $\R^n$ be split into at most $n+1$ parts
of smaller diameter? Borsuk \cite{Borsuk} asked this question in 1933.
The vertices of a unit equilateral triangle need three parts, since any
two are a unit apart. For the same reason, the vertices of a regular
$n$-simplex need $n+1$ parts.

Write $b_n$ for the least number of parts that always suffices. Borsuk proved
$b_2=3$. Perkal \cite{Perkal} announced $b_3=4$ in 1947, and Eggleston
\cite{Eggleston} published a proof in 1955. In 1993, Kahn and Kalai
\cite{KK} disproved the conjecture in high dimensions. Four is the first
dimension in which the question remains open.

Lassak's 1982 bound $b_n\le2^{n-1}+1$ gives $b_4\le9$ \cite{Lassak}.
In May 2026, Tolmachev and Voronov \cite[Theorem~2]{TV} reduced this
to eight. In October 2026, Arman, Bondarenko, Prymak and Radchenko
\cite[Theorem~1]{ABPR} covered every set of diameter one in $\R^4$
with seven balls of diameter less than $0.9999$. We prove the following bound.

\begin{theorem}\label{thm}
Every bounded subset of $\R^4$ of positive diameter admits a partition into
at most six sets of strictly smaller diameter. In particular, $5\le b_4\le6$.
\end{theorem}

Theorem~\ref{thm} is formalised in Lean~\cite{Lean}.

A part of diameter less than one need not lie in a ball of diameter
less than one. We use this freedom to improve the seven-ball construction
of Arman, Bondarenko, Prymak and Radchenko. We retain their division
according to circumradius and their reduction to three contact points
at intermediate radii.

The main change is at large circumradius. Their five balls centred at
rescaled contacts leave two regions covered by two further balls.
We instead choose five balls so that no two uncovered points can be
at distance one. The entire complement is then a single sixth part.
The second-moment estimate of \cite[Lemma~6]{ABPR}, applied to both
endpoints of a hypothetical diametral pair, gives the contradiction.

At small circumradius, six cones partition the enclosing ball.
At intermediate radii, we divide the outer sets determined by contact
triangles using six affine functions. Enclosing each region in a polytope
reduces its diameter bound to distances between vertices. The choice of
contact triangle also matters: some triangles satisfying the initial
parameter bounds cannot arise from successive nearest facets.
Lemma~\ref{select} excludes these cases.

By scaling and taking closures it suffices to partition a compact set $K$ of
diameter one. Centre its smallest enclosing ball at zero and write $r$ for its
squared radius. There are affinely independent contacts $p_1,\ldots,p_m\in K$
and weights $\lambda_i>0$ with $m\le5$, $\norm{p_i}^2=r$,
$\sum_i\lambda_i=1$ and $\sum_i\lambda_i p_i=0$ \cite[Lemma~2]{ABPR}. Thus
\begin{equation}\label{jung}
 r=\sum_{i<j}\lambda_i\lambda_j\norm{p_i-p_j}^2
 \le \tfrac12\bigl(1-\textstyle\sum_i\lambda_i^2\bigr)
 \le \frac{m-1}{2m}\le\frac25.
\end{equation}
Put $r_-=153/500$ and $r_+=1957/5000$. Let $B(c,t)$ be the open ball of radius $t$ about $c$.

\section{Large circumradius}
We seek five balls whose compact complement in $K$ has no pair at distance
one. Avoiding the balls bounds the combined second moment of such a pair
from above; its distance gives a contradictory lower bound.

\begin{proposition}\label{large}
Suppose $r\ge r_+$. There are five contacts $p_i$ such that
\[
 F=K\setminus\bigcup_{i=1}^5 B(\alpha p_i,\rho),\qquad
 \alpha=\tfrac25,\quad \rho=\tfrac{9999}{20000},
\]
contains no pair at distance one.
\end{proposition}

\begin{proof}
Since $r>3/8$, \eqref{jung} gives $m=5$. Put
\[
 M=\sum_i\lambda_i p_i p_i^{\mathsf T},\quad
 \Phi(x)=x^{\mathsf T}Mx,\quad U=\tfrac12-r,\quad N=UI-M.
\]
By \cite[Lemma~6]{ABPR}, $M\preceq UI$. Hence $N$ is positive semidefinite
and $\tr N=4U-r=2-5r$.

Let $x\in F$, $q=\norm{x}^2$, and $\xi=\alpha^2r-\rho^2$.
The inequalities $\norm{x-p_i}\le1$ and $\norm{x-\alpha p_i}\ge\rho$ give
\[
 -\frac{1-r-q}{2}\le\ip{x}{p_i}\le\frac{q+\xi}{2\alpha}.
\]
The scalar products have weighted mean zero. For $a\le s\le b$ we have
$(s-a)(b-s)\ge0$. Averaging this inequality gives
\begin{equation}\label{one}
 \Phi(x)\le\frac{(1-r-q)(q+\xi)}{4\alpha}.
\end{equation}

Suppose now that $x,y\in F$ satisfy $\norm{x-y}=1$. Put $t=\norm{y}^2$,
$Q=q+t$ and $D=q-t$. Then $\ip{x}{y}=(Q-1)/2$ and $Q\le2r<1$.
Adding \eqref{one} for $x$ and $y$ gives
\[
 \Phi(x)+\Phi(y)\le V(Q,D):=
 \frac{(1-r)(Q+2\xi)-(Q^2+D^2)/2-\xi Q}{4\alpha}.
\]
For a lower bound, let $S=xx^{\mathsf T}+yy^{\mathsf T}$. The Gram matrix
of $x,y$ gives
\[
 \lambda_{\max}(S)=\tfrac12\bigl(Q+\sqrt{D^2+(1-Q)^2}\bigr).
\]
Write $\kappa=(2-5r)/2$. Since $\tr(NS)\le2\kappa\lambda_{\max}(S)$,
\[
 \Phi(x)+\Phi(y)=UQ-\tr(NS)
 \ge (U-\kappa)Q-\kappa\sqrt{D^2+(1-Q)^2}.
\]
Subtract $V(Q,D)$ from this lower bound and use
$\sqrt{a^2+D^2}\le a+D^2/(2a)$ with $a=1-Q$. The difference is at least
\[
 G(Q,r)+D^2\left(\frac1{8\alpha}-\frac{\kappa}{2(1-Q)}\right),
 \qquad G(Q,r)=UQ-\kappa-V(Q,0).
\]
Since $r\ge3/10$, we have $4\alpha\kappa\le1-2r\le1-Q$, so the
coefficient of $D^2$ is nonnegative.

It remains to show $G(Q,r)>0$. Completing the square gives
\[
 8\alpha G(Q,r)=(Q+B/2)^2+\Delta(r),\qquad \Delta=C-B^2/4,
\]
where $B=8\alpha U-2(1-r-\xi)$ and $C=-8\alpha\kappa-4(1-r)\xi$.
On $[r_+,2/5]$ we have $-2<B<0$, $B'=-22/25$ and $C'>6$, so
$\Delta'>5$. Substitution gives $\Delta(r_+)>1/10000$.
Thus $G(Q,r)>0$, contradicting the upper bound for $\Phi(x)+\Phi(y)$.

A nonempty $F$ is compact and contains no unit pair, so $\diam F<1$.
Assigning each point of $K\setminus F$ to the first ball containing it
gives five parts of diameter at most $2\rho<1$.
\end{proof}

\section{Small circumradius}
\begin{figure}[htbp]
\centering
\includegraphics[width=.94\linewidth]{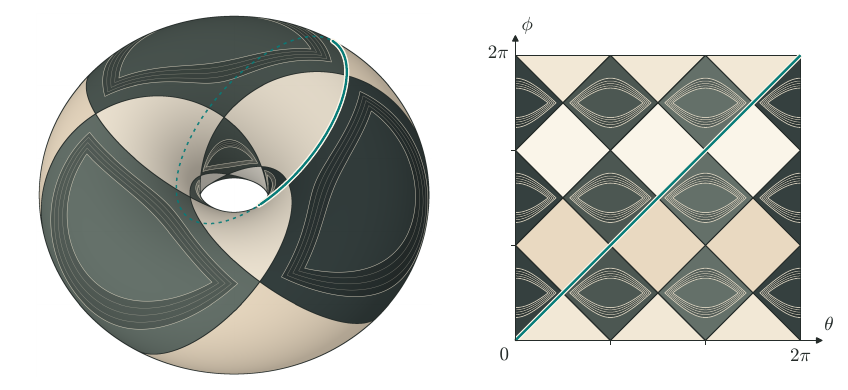}
\caption{Stereographic projection of the Clifford torus (left) and its
angle square (right). The boundaries $\theta\pm\phi\in(2\pi/3)\mathbb Z$
project to six Villarceau circles. The turquoise diagonal $\theta=\phi$ traces
one such circle, dashed where hidden. In the dark regions $m(\theta)>m(\phi)$.
The fine contours mark $m(\theta)-m(\phi)=k/32$, $4\le k\le8$.}
\label{fig:torus}
\end{figure}

For small circumradius we partition the entire circumball using two
equilateral triples in orthogonal planes. Let $d_1,\ldots,d_6$ be the unit vectors at angles $0,2\pi/3,4\pi/3$
in the two coordinate planes of $\R^2\oplus\R^2$. The six cones on which
the respective $\ip{x}{d_i}$ are maximal are congruent.
On the Clifford torus in Figure~\ref{fig:torus} the two planar components
are $2^{-1/2}(\cos t,\sin t)$, with respective angles $t=\theta,\phi$.
Their maximal scalar products are $m(\theta)/\sqrt2$ and $m(\phi)/\sqrt2$,
where $m(t)=\max_{k=0,1,2}\cos(t-2\pi k/3)$.

\begin{proposition}\label{small}
A ball of radius $R$ in $\R^4$ has a partition into at most six sets of squared
diameter at most $(2+2\sqrt{2/5})R^2$. This bound is less than one for
$R^2\le r_-$.
\end{proposition}

\begin{proof}
In the cone for $d_1=(1,0,0,0)$ every ray has a representative
$(1,b,z)$, where
\[
 |b|\le\sqrt3,\qquad z\in T=\conv\{(-2,0),(1,\sqrt3),(1,-\sqrt3)\}.
\]
Indeed, the inequalities $\ip{x}{d_1}\ge\ip{x}{d_j}$ force $x_1>0$
on every nonzero ray and give this section at $x_1=1$.

Fix $z,z'\in T$ and put
$A=1+\ip{z}{z'}$, $B=1+\norm{z}^2$, $B'=1+\norm{z'}^2$.
For fixed $s=|b|$ and $t=|b'|$ the smaller cosine is
\[
 f(s,t)=\frac{A-st}{\sqrt{(B+s^2)(B'+t^2)}}.
\]
If $s\le t$, the sign of $\partial_s f$ is the sign of
$-B(t-s)-(B+A)s$. This is nonpositive since $\norm{z'}\le2$ and
$B+A=2+\norm{z+z'/2}^2-\norm{z'}^2/4\ge1$.
By symmetry we can increase the smaller coordinate to the larger, then
increase both to $\sqrt3$, without increasing $f$. The resulting cosine is
\[
 f(\sqrt3,\sqrt3)=
 \frac{\ip{z}{z'}-2}{\sqrt{(4+\norm{z}^2)(4+\norm{z'}^2)}}.
\]

At a minimum with $z$ interior to $T$, differentiation gives
$(4+\norm{z}^2)z'=(\ip{z}{z'}-2)z$, whence $z=-2z'$.
Then $v=\norm{z'}^2\le1$ and the cosine is
$-\sqrt{(1+v)/(4+v)}\ge-\sqrt{2/5}$.
Otherwise symmetry puts $z=(1,t)$ on an edge. Writing $z'=(a,b)$
gives $3b^2\le(a+2)^2$. By Cauchy--Schwarz,
\[
 (a-2+bt)^2\le (5+t^2)\bigl((a-2)^2/5+b^2\bigr)
 \le \tfrac25(5+t^2)(4+a^2+b^2).
\]
This gives the same lower bound. Consequently two points in one cone
of the ball are at squared distance at most $(2+2\sqrt{2/5})R^2$.
Assign ties to the least index. Since $\sqrt{2/5}<19/30$, at $R^2=r_-$
the bound is less than $(49/15)(153/500)=2499/2500$.
\end{proof}

The rays through $(1,\sqrt3,1,0)$ and $(1,-\sqrt3,-2,0)$ attain the
angular bound, so the constant is sharp for this partition.

\section{Intermediate circumradius}
Suppose $r_-\le r\le r_+$. Since $r_->1/4$, \eqref{jung} gives at least
three contacts. Take successive facets nearest the circumcentre of the
current simplex until a triangle remains. By \cite[Lemma~3 and Proposition~4]{ABPR}, the
triangle is acute, and its squared circumradius $\tau$ and the distance
$h$ of its plane from zero satisfy
\[
 5r/6\le\tau\le\min(r,1/3),\qquad h^2=r-\tau.
\]
Indeed, passing from a $k$-simplex to its nearest facet retains at least
$1-1/k^2$ of its squared circumradius. The factor $5/6$ is the product
$(15/16)(8/9)$ when both steps are needed. For the selected contacts
$P_1,P_2,P_3$, we have
\begin{equation}\label{H}
 K\subset H=\clB(0,\sqrt r)\cap\bigcap_{i=1}^3\clB(P_i,1).
\end{equation}

\subsection{Contact triangles}
Order the angles $A\ge B\ge C$ and write $a=\cos A$, $b=\cos B$.
The longest side is at most one, and $(\pi-A)/2\le B\le A$. Thus
\[
 a_0:=\sqrt{\max(0,1-1/(4\tau))}\le a\le\tfrac12,\qquad
 a\le b\le\sqrt{(1-a)/2}.
\]
Following \cite[Lemma~10]{ABPR}, for $u\in[0,1]^4$ put
\begin{align*}
 t&=u_1\sqrt{r_+-1/3}-(1-u_1)\sqrt{1/3-r_-}, & r&=1/3+t|t|,\\
 h&=(1-u_2)\sqrt{r/6}+u_2\max(t,0), & \tau&=r-h^2,\\
 a&=(1-u_3)/2+u_3a_0, & b&=(1-u_4)\sqrt{(1-a)/2}+u_4 a.
\end{align*}
Each interpolation covers its closed interval. Put
$w_\pm(c)=(2c^2-1,\pm2c\sqrt{1-c^2})$. After an orthogonal change of
coordinates the contacts are
\[
 (\sqrt\tau,0,h,0),\qquad
 (\sqrt\tau\,w_+(a),h,0),\qquad
 (\sqrt\tau\,w_-(b),h,0).
\]
The side lengths are $2\sqrt\tau$ times the sines of $A,B,C$, and the
plane has distance $h$ from zero. Write $H_u$ for \eqref{H} in these
coordinates. Rotating the last two axes carries $H_u$ and its partition
to congruent sets.

A selected triangle must also arise from successive nearest facets.
The following restriction on its shape will exclude the boxes left over
by the polyhedral partitions.

\begin{lemma}\label{select}
Suppose the selected triangle has $T=\tau/r<8/9$. For every $q$ such
that $0<q\le1/3$ and $q^2\le1-16T/15$, we have
\begin{equation}\label{sel}
 a^2\ge\frac{7-9U_q-9E_q}{16},\qquad
 U_q=\frac{1-2q}{3q},\qquad E_q=\sqrt{2-\frac{6q^2}{1-2q}}.
\end{equation}
\end{lemma}

\begin{proof}
With three or four initial contacts we would have $\tau\ge8r/9$.
Hence the triangle is a nearest facet of a tetrahedron whose squared
circumradius $S$ satisfies $15r/16\le S\le9\tau/8$. Every facet of
this tetrahedron is at distance at least $\sqrt{S-\tau}\ge q\sqrt S$
from its circumcentre.

Translate that centre to zero. Write $v_i$ for the vertices, $n_i$ for
the outward unit facet normals, $d_i$ for the corresponding distances,
and $\lambda_i$ for the barycentric weights of zero. The facet opposite
$v_i$ gives
\[
 -\ip{n_i}{v_i}=\frac{1-\lambda_i}{\lambda_i}d_i\le\sqrt S .
\]
Thus $\lambda_i\ge q/(1+q)$, and summing over the other three weights
gives $\lambda_i\le(1-2q)/(1+q)$. These two bounds imply
\[
 d_i/\sqrt S\le U_q,\qquad
 -\ip{n_i}{v_i/\sqrt S}\ge 3q^2/(1-2q),\qquad
 \norm{v_i/\sqrt S+n_i}\le E_q.
\]
For $j\ne i$ we have $\ip{n_i}{v_j}=d_i$, so
$\ip{v_i}{v_j}/S\ge-U_q-E_q$.
In particular the longest side of the triangle satisfies
\[
 4\tau(1-a^2)\le S(2+2U_q+2E_q)
                  \le(9\tau/8)(2+2U_q+2E_q).
\]
Rearrangement gives \eqref{sel}. The argument allows ties between nearest facets.
\end{proof}

\subsection{Polyhedral partitions}
On a closed dyadic box $Q\subset[0,1]^4$, let
$s_j(x)=\ip{c_j}{x}+\beta_j$, $1\le j\le6$, have rational coefficients.
Assign $x\in H_u$ to the least $j$ maximizing $s_j(x)$. For every $n$,
\eqref{H} gives
\begin{equation}\label{support}
 \ip{n}{x}\le\sqrt r\,\norm n,\qquad
 \ip{n}{x}\le\ip{n}{P_i}+\norm n\quad(i=1,2,3).
\end{equation}
For rational $n$, a rational upper bound on any right side throughout $Q$
gives a halfspace containing every $H_u$, $u\in Q$. Intersect finitely
many such halfspaces with $[-1,1]^4$ and the five halfspaces $s_k\le s_j$.
The resulting rational polytope $C_j$ contains the $j$th region. Its
diameter is the largest distance between its vertices, since a difference
of convex combinations of vertices is a convex combination of their
pairwise differences.

\begin{proposition}\label{middle}
The cube $[0,1]^4$ admits a finite cover by closed dyadic boxes.
On each box, either six affine functions partition every $H_u$ into
sets of squared diameter at most $1-86\cdot10^{-10}$, or no $u$
represents a triangle selected by successive nearest facets.
\end{proposition}

\begin{proof}
The supplementary data list the four intervals defining
each box. For each of the $1512$ partition entries they give a bound $M_Q$
and the source of the six affine functions and their supporting halfspaces.
To bound \eqref{support}, we use rational intervals with outward rounding
and intersect the bounds from the following identities and extrema.
Write $u=u_2$, $A=\sqrt{r/6}$, $B=\sqrt{\max(r-1/3,0)}$ and $d=A-B$.
Then $h=A-ud=B+(1-u)d$, and
\[
 \tau=\frac{5r}{6}+2uAd-u^2d^2
      =\min(r,1/3)-2(1-u)Bd-(1-u)^2d^2.
\]
The function $\tau(r,u)$ increases in $u$. For fixed $u$, it is linear
in $r\le1/3$ and convex in $r\ge1/3$. For $0<u<1$, the latter branch
has its minimum at $r_*(u)$, where $p(u)=25+70u+25u^2$ and
\begin{align*}
 r_*(u)&=\frac{1+(5+7u)/\sqrt{p(u)}}6,\\
 \tau(r_*(u),u)&=\frac{u^2}3+
 \frac{(1-u)(5+7u+\sqrt{p(u)})}{36}.
\end{align*}
Thus the extrema on a rectangle occur at its $u$-endpoints and the
$r$-endpoints, at $r=1/3$, or at $r_*(u)$ when included. For $u=0,1$
the formulas reduce to $5r/6$ and $\min(r,1/3)$. These bounds are
intersected with $5r/6\le\tau\le\min(r,1/3)$.

The parameter $a$ increases in $\tau$ and decreases in $u_3$.
For fixed $w=u_4$, the function
$b=(1-w)\sqrt{(1-a)/2}+wa$ is concave in $a$ and decreases in $w$.
Its extrema are bounded using endpoints and, for $w>0$, the stationary
point $a=1-(1-w)^2/(8w^2)$ when included; the value there is
$w+(1-w)^2/(8w)$.

We also retain the dependence of the first contact coordinates on $\tau$.
Put $s=\sqrt\tau$, $v=u_3$ and
\[
 a(s,v)=\frac{1-v}{2}+v\sqrt{1-\frac1{4s^2}},\qquad
 X(s,v)=s(2a(s,v)^2-1).
\]
On $1/2\le s\le1/\sqrt3$, $X$ is concave in $s$ and decreases in $v$.
Its minimum on a rectangle is therefore at an $s$-endpoint with maximal
$v$. Moreover, $\partial_sX(1/\sqrt3,v)=3v-1/2$, so $X$ increases in
$s$ if $v\ge1/6$; for $v=0$ it equals $-s/2$.
These give upper bounds in those cases; otherwise we use interval evaluation.
The first coordinates satisfy
$(P_2)_1=X\le(P_3)_1\le-sa$. The upper bound for $-sa$ is attained at the
smallest $s$ and largest $v$. All these assertions follow by differentiation.

For each normal we take the smaller of the bound from the separate
coordinate intervals and the following bound. For $P_2$ and $P_3$ the
planar contribution to $\ip{n}{P_i}$ is
$\sqrt\tau$ times
\[
 n_1(2z-1)\pm 2n_2\sqrt{z(1-z)},\qquad z=a^2\ \text{or}\ b^2.
\]
The endpoints and any interior maximum give an upper bound $C$ for this
factor; for $P_1$ take $C=n_1$. It remains to bound
$C\sqrt{r-h^2}+n_3h$. When $h_{\max}^2<r_{\min}$,
use $r_{\max}$ for $C\ge0$ and $r_{\min}$ for $C<0$, then maximize in $h$.
Otherwise use separate bounds for $\sqrt\tau$ and $h$.

Cutting $[-1,1]^4$ by the halfspaces gives the vertices of $C_j$. Put
$L_i(v)=\lfloor2^{25}v_i\rfloor$. Since
$2^{25}|v_i-w_i|\le |L_i(v)-L_i(w)|+1$, a common squared-diameter bound is
\[
 M_Q=2^{-50}\max_{j,v,w}
          \sum_{i=1}^4\bigl(\lvert L_i(v)-L_i(w)\rvert+1\bigr)^2 .
\]
Here $v,w$ range over the vertices of the same $C_j$.
For a partition obtained by $x=z+\lambda y$ from one with bound $M_*$,
take $M_Q=\lambda^2M_*$. For each retained normal $n$, the support bounds
satisfy $\widehat b_n\le\ip{n}{z}+\lambda b_n$, so $\ip{n}{y}\le b_n$.
Also $\sqrt{r_{\max}}+\norm z_1\le\lambda$, giving $y\in[-1,1]^4$.

The box $[15/16,1]\times[1/8,1/4]\times[3/4,1]\times[1/2,1]$ gives the largest $M_Q$.
Its reference partition and dilation give
\[
 M_*=\frac{1069228342837886}{2^{50}},\qquad
 \lambda=\frac{1026159}{10^6},\qquad
 \lfloor10^{10}(1-\lambda^2M_*)\rfloor=86.
\]

Subdivide each of the ten slabs perpendicular to the $u_1$-axis at
the box endpoints in each coordinate. A closed cell $\prod_i[l_i,h_i]$
lies in $\prod_i[a_i,b_i]$ precisely when $a_i\le l_i\le h_i\le b_i$
for all $i$. For each interval $I$ below, $P$ counts cells in partition
boxes and $E$ counts the remaining cells. The last column is
$\lfloor10^{10}(1-\max_Q M_Q)\rfloor$ over the partition entries
in those slabs.
\begin{center}
\begin{tabular}{@{}lrrrr@{}}
\hline
$I$ for $u_1$ & All cells & $P$ & $E$ & $\lfloor10^{10}(1-\max_Q M_Q)\rfloor$\\
\hline
$[0,11/16]$    & $39\,248$ & $39\,248$ & 0  & $12\,008$\\
$[11/16,3/4]$  & $3\,150$  & $3\,138$  & 12 & $48\,119$\\
$[3/4,1]$      & 265       & 265       & 0  & 86\\
\hline
\end{tabular}
\end{center}
The twelve remaining cells lie in the union of nine boxes on each of
which $T<8/9$. On each box a rational $q$ satisfies the hypotheses of
Lemma~\ref{select}, but the upper bound for $a^2$ is strictly below
the right side of \eqref{sel}. Hence these boxes contain no selected
triangle. The last column gives $M_Q\le1-86\cdot10^{-10}$ throughout
the partition boxes.
\end{proof}

\begin{proof}[Proof of Theorem~\ref{thm}]
Propositions~\ref{large} and~\ref{small} give six parts when
$r\ge r_+$ or $r\le r_-$. Otherwise select a contact triangle.
Its parameter $u$ lies in a partition box by Proposition~\ref{middle}.
The six regions of $H_u$ restrict to $K$ by \eqref{H}. Restricting to
the original set and rescaling proves the upper bound. The regular
four-simplex gives the lower bound.
\end{proof}

\subsection{Construction of the partitions}
We alternate between adjusting the partition and sharpening its polyhedral
diameter bounds. The affine coefficients determine the regions; additional
supporting halfspaces cut away unwanted vertices of the enclosing polytopes.
For fixed supporting halfspaces, let $d_j$ be the largest squared
vertex distance in $C_j$. We vary the affine coefficients to reduce
$\max_j d_j$, using L-BFGS-B with the approximation
\[
 \frac1\eta\log\left(\sum_{j=1}^6 e^{\eta d_j}\right)
\]
and increasing $\eta$ in stages. Subtracting the mean coefficient row
removes the freedom to add a common affine function; a penalty on the
coefficient norm controls their common scale. At each evaluation, linear
programming supplies interior points for halfspace intersection. A farthest
vertex pair determines $d_j$; differentiating four independent active
facet equations at each endpoint gives the derivatives used in the search.
The best coefficients are retained according to $\max_j d_j$.

Initial supporting normals are obtained by rounding vectors from the
ball centres to sampled boundary points. Partitions on neighbouring
boxes supply initial affine coefficients, and subdivision tightens the
support bounds. Some partitions are obtained by translation and dilation,
using the inequalities in the proof of Proposition~\ref{middle}.

For each region we retain supporting halfspaces incident to its vertices.
To refine a partition, we examine a pair of vertices at the largest distance.
For an endpoint $v$, rational approximations to $v$ and $v-P_i$, with
contacts evaluated at corners of the parameter box, give candidate normals.
We retain inequalities violated by $v$, with support bounds valid
throughout the box. This cuts away unwanted vertices without removing
any point of $H_u$. The affine coefficients can then be varied again.
In these continuations the coefficients are rounded to multiples of
$2^{-24}$, after which the incident halfspaces are selected again.
The diameter bounds in Proposition~\ref{middle} are computed from the
resulting rational polytopes.

The $1512$ partition entries use $1501$ stored
partitions. The table counts each entry once, including reused partitions.
Halfspaces shared by several regions are counted once per box in the
first row and once per region in the second.
\begin{center}
\begin{tabular}{@{}lrrr@{}}
\hline
Supporting halfspaces & Minimum & Mean & Maximum\\
\hline
Per parameter box & 188 & 285.90 & 4096\\
Per region & 31 & 69.68 & 1712\\
\hline
\end{tabular}
\end{center}
Each region also has eight cube inequalities and five affine comparison
inequalities: $44$--$1725$ inequalities in total, with mean $82.68$.

\section{Concluding remarks}
If $b_4=6$, four is the least dimension in which Borsuk's conjecture fails.

\smallskip\noindent\emph{Can a compact set of diameter one on a sphere
of radius $R\le1/\sqrt2$ in $\R^4$ require six parts of smaller diameter?}

\smallskip
Such a set would give $b_d\ge d+2$ for every $d\ge4$ by the construction
of Hinrichs and Richter \cite[Lemma~9]{HR}, also used for the
nine-dimensional example in \cite[Corollary~7.1]{OAI}.
Suppose $K$ is such a set and centre its sphere at zero. To
$K\times\{0\}$ adjoin $(0,h)$ in one
more dimension, where $h=\sqrt{1-R^2}$. This point is at distance one
from every point of $K\times\{0\}$, so it requires a new part. The
enlarged set lies on the sphere with centre $(0,h-1/(2h))$ and radius
$1/(2h)\le1/\sqrt2$. Iterating would leave $d=1,2,3$ as precisely the
dimensions in which $d+1$ parts always suffice.

\section*{Acknowledgements}
The author thanks Andriy Prymak and Vsevolod Voronov for helpful comments on an earlier draft.

\section*{AI Declaration}
OpenAI's Astra assisted with exploring proofs, preparing the Lean formalisation, and
in the article preparation process. The author takes full responsibility for the article.

\section*{Mathematical software}
Python programs were used to construct the partitions and to compute the
rational bounds in Proposition~\ref{middle}. The full theorem was checked
in Lean~4.36.0-rc2 with mathlib; the formalisation is cited in~\cite{Lean}.

\section*{Data availability}
The data and Lean formalisation are available at \href{https://doi.org/10.5281/zenodo.23966881}{doi:10.5281/zenodo.23966881}.

\section*{Conflicts of Interest}
The author declares no conflicts of interest.


\begin{thebibliography}{10}
\bibitem{ABPR}
A.~Arman, A.~Bondarenko, A.~Prymak and D.~Radchenko,
\emph{Borsuk number in four dimensions is at most seven},
\href{https://arxiv.org/abs/2610.39078v1}{arXiv:2610.39078} (2026).
\bibitem{Lean}
L.~R.~Chalmers, \emph{Lean formalisation of the four-dimensional six-part theorem}
(2026), \href{https://doi.org/10.5281/zenodo.23966881}{doi:10.5281/zenodo.23966881}.
\bibitem{Borsuk}
K.~Borsuk, \emph{Drei S\"atze \"uber die $n$-dimensionale euklidische
Sph\"are}, Fund. Math. \textbf{20} (1933), 177--190.
\bibitem{Eggleston}
H.~G.~Eggleston, \emph{Covering a three-dimensional set with sets of
smaller diameter}, J. London Math. Soc. \textbf{30} (1955), 11--24.
\bibitem{HR}
A.~Hinrichs and C.~Richter, \emph{New sets with large Borsuk numbers},
Discrete Math. \textbf{270} (2003), 137--147.
\bibitem{KK}
J.~Kahn and G.~Kalai, \emph{A counterexample to Borsuk's conjecture},
Bull. Amer. Math. Soc. (N.S.) \textbf{29} (1993), 60--62.
\bibitem{Lassak}
M.~Lassak, \emph{An estimate concerning Borsuk partition problem},
Bull. Acad. Polon. Sci. S\'er. Sci. Math. \textbf{30} (1982), 449--451.
\bibitem{OAI}
OpenAI, \emph{A nine-dimensional counterexample to Borsuk's covering
assertion}, \href{https://github.com/openai/math/blob/main/preprints/A-nine-dimensional-counterexample-to-Borsuks-covering-assertion-September-23-2026/paper.pdf}{preprint} (2026).
\bibitem{Perkal}
J.~Perkal, \emph{Sur la subdivision des ensembles en parties de
diam\`etre inf\'erieur}, announcement in the Wroc{\l}aw proceedings,
\href{https://doi.org/10.4064/cm-1-1-37-50}{Colloq. Math. \textbf{1}
(1947), 45}.
\bibitem{TV}
A.~Tolmachev and V.~Voronov, \emph{Reducing the upper bound for the
Borsuk number in $\R^4$ to 8},
\href{https://arxiv.org/abs/2605.19068v1}{arXiv:2605.19068} (2026).
\end{thebibliography}
\end{document}